# Continuous families of isospectral metrics on simply connected manifolds

By Dorothee Schueth*

## Abstract

We construct continuous families of Riemannian metrics on certain simply connected manifolds with the property that the resulting Riemannian manifolds are pairwise isospectral for the Laplace operator acting on functions. These are the first examples of simply connected Riemannian manifolds without boundary which are isospectral, but not isometric. For example, we construct continuous isospectral families of metrics on the product of spheres $S^4 \times S^3 \times S^3$. The metrics considered are not locally homogeneous. For a big class of such families, the set of critical values of the scalar curvature function changes during the deformation. Moreover, the manifolds are in general not isospectral for the Laplace operator acting on 1-forms.

## Introduction

Spectral geometry deals with the mutual influences between the geometry of a Riemannian manifold and the spectrum of the associated Laplace operator acting on smooth functions. Until 1964 it was not known whether the spectrum determines the geometry completely. Then J. Milnor constructed the first counterexample, namely, a pair of isospectral, nonisometric flat tori in dimension sixteen. Numerous other examples of isospectral manifolds followed in the 1980's and 1990's. See, e.g., [Vi], [Ik], [GW2], [Su], [BT], [Gt1] for examples of discrete families of isospectral manifolds, or [GW1], [DG1,2], [BG], [Sch], [Gt2] for isospectral deformations, that is, continuous families of isospectral manifolds.

In all the examples cited above, the isospectral manifolds arise as different quotients (by discrete groups of isometries) of a common Riemannian covering

1991 *Mathematics Subject Classification.* 58G25. *Key words and phrases.* Laplacian, isospectral deformations, simply connected manifolds, products of spheres, scalar curvature, 1-form spectrum, heat invariants.

*The author is partially supported by SFB 256, Bonn. Moreover, the author wishes to thank the GADGET program for supporting her participation in the NSF-CNRS workshop on Spectral Geometry in Grenoble (1997) where part of this work was carried out.

manifold. In particular, such quotients can neither differ in their local geometry nor be simply connected.

In 1993, C. Gordon [Go1], [Go2] constructed the first pairs of isospectral closed Riemannian manifolds which are not locally isometric. S. Szabó [Sz] had constructed the first such examples in the case of manifolds with boundary even earlier, but did not publish them at that time. In 1996, even continuous families of isospectral metrics with different local geometries were obtained by C. Gordon and E. Wilson [GW3] in the case of manifolds with boundary. More precisely, in their examples the underlying manifolds are products of a ball with a torus. Recently it was discovered [GGSWW] that also the boundaries of these manifolds are isospectral, which led to the first examples of isospectral deformations with changing local geometry in the case of closed manifolds, namely, on products of a sphere with a torus. Szabó [Sz] independently gave finite families of isospectral manifolds which represent special cases of the examples in [GW3] and [GGSWW].

However, there were still no examples known, not even pairs, of *simply connected*, closed isospectral manifolds. Simply connected examples existed until now only in the case of manifolds with (nonsmooth) boundary, most notably the famous pairs of isospectral planar domains by Gordon, Webb, and Wolpert [GWW].

In the present paper, we construct the first examples of simply connected, closed Riemannian manifolds which are isospectral, but not isometric; we even give continuous families of such manifolds. We will show:

MAIN THEOREM. *Let $S^n$ be a sphere of dimension $n \geq 4$, and let $S$ be a compact, simply connected Lie group whose maximal torus is of dimension at least two (for example, $S := \mathrm{SU}(2) \times \mathrm{SU}(2) \cong S^3 \times S^3$). Then there exist continuous isospectral families of Riemannian metrics on $S^n \times S$ which are not locally isometric; more precisely, the set of critical values of the corresponding scalar curvature functions changes during the deformations.*

In particular, we obtain such families of metrics on $S^n \times S^3 \times S^3$ for every $n \geq 4$. These isospectral deformations can take place arbitrarily close to the standard metric on $S^n \times S^3 \times S^3$ (although none of our deformations is one of the standard metric itself).

Our construction is related to, and motivated by, the isospectral examples in [GGSWW], where the underlying manifolds are products of a sphere with a torus. The idea is to embed the nonsimply connected factor, that is, the torus, into a compact, simply connected Lie group and to extend the metrics to the bigger manifolds in such a way that they remain isospectral.

Concerning the change in geometry during the deformations, it turns out that in many cases, the isospectral manifolds can be distinguished by the set of critical values assumed by the scalar curvature function. Note that this function is in general nonconstant here since the manifolds are not locally homogeneous.

Moreover, our isospectral deformations are in general not isospectral for the Laplace operator acting on 1-forms. More precisely, the manifolds can be distinguished already by the heat invariants of the Laplacian acting on 1-forms. This is interesting in several respects: There are only few examples known of manifolds which are isospectral on functions, but not on 1-forms. The only known examples are those given in [Ik] (lens spaces), [Go3] (Heisenberg manifolds), and [Gt1], [Gt2] (certain three-step nilmanifolds). Among these, only Gornet's examples include continuous families of such manifolds [Gt2]. In all these examples, the isospectral manifolds arise as different quotients of some common Riemannian covering. In particular, these isospectral manifolds, even though they are not isospectral on 1-forms, always share the same heat invariants of the Laplacian acting on 1-forms. The isospectral families constructed in the present paper are the first ones where this is not the case.

We will construct a large class of continuous isospectral families in which the manifolds can be distinguished by either of the two criteria mentioned above (critical values of the scalar curvature, resp. spectrum on 1-forms).

The paper is organized as follows:

In Section one, we introduce the Laplace operator and a general isospectrality principle (Proposition 1.2) which was discovered by C. Gordon [Go2] and is the key of the isospectrality proof in Section two, as it was the case also for the manifolds constructed in [GW3] and [GGSWW]. We then prepare our construction of isospectral manifolds by introducing a certain class of metrics $g_j$ on $\mathbb{R}^{n+1} \times S \supset S^n \times S$, where $S$ is a compact, simply connected Lie group, and discussing some of their properties. These metrics are related to left invariant metrics on two-step nilpotent Lie groups diffeomorphic to $\mathbb{R}^{n+1} \times H$, where $H$ is a closed abelian subgroup of $S$. Each metric $g_j$ is associated with a linear map $j : \mathbb{R}^r \to so(n+1)$, where $r$ is the dimension of $H$.

In Section two, we consider the submanifolds $S^n \times S$ of $\mathbb{R}^{n+1} \times S$, endowed with the restricted metrics, and show that under a certain condition on $j$ and $j'$ which was introduced in [GW3], the Riemannian manifolds $(S^n \times S, g_j)$ and $(S^n \times S, g_{j'})$ are isospectral (Theorem 2.2). We finish this section by giving a class of examples of continuous families of isospectral metrics in the case $n = 4$.

In Section three, we study the geometry of the isospectral manifolds constructed in Section two. We establish a formula for the scalar curvature function on these manifolds (Lemma 3.3). This allows us to describe the sets of critical points and critical values of the scalar curvature in terms of eigenvectors and eigenvalues of the Ricci tensor of the corresponding nilpotent Lie group which we mentioned above (Corollary 3.4) and thus to obtain a suffi-cient condition for nonisometry of two isospectral manifolds $(S^n \times S, g_j)$ and $(S^n \times S, g_{j'})$.

In Section four, we show that for $n \geq 4$, any generic linear map $j : \mathbb{R}^2 \to \mathrm{so}(n+1)$ is contained in a continuous family $j(t)$ of such maps with the property that the corresponding manifolds $(S^n \times S, g_{j(t)})$ are isospectral and pairwise satisfy the nonisometry condition established in Section three. More precisely, the manifolds have pairwise different sets of critical values of the scalar curvature. The families $j(t)$ are obtained as the flow lines of a certain vector field on the space of all linear maps from $\mathbb{R}^2$ to $\mathrm{so}(n+1)$.

Section five is devoted to the spectrum on 1-forms. Using the heat invariants for the Laplace operator on 1-forms and on functions, we show that in each of the special families constructed in Section four, the spectrum on 1-forms changes during the deformation.

The author wishes to thank Carolyn Gordon and Ruth Gornet for interesting discussions. Concerning the 1-form spectrum, she also wishes to acknowledge that Ruth Gornet was the originator of the idea to study the heat invariants for the form spectra of the manifolds constructed in [GGSWW], thereby also inspiring Section five of the present paper.

## 1. Preliminaries

1.1 *Definition.* Let $(M, g)$ be a compact Riemannian manifold, and let $\Delta_g$ be the Laplacian acting on functions by

$$(\Delta_g f)(p) := -\sum_{i=1}^n \frac{d^2}{dt^2}\Big|_{t=0} f(c_i(t)) \qquad \text{for } p \in M,$$

where the $c_i$ are geodesics starting in $p$ such that $\{\dot{c}_1(0), \ldots, \dot{c}_n(0)\}$ is an orthonormal basis for $T_p M$. The discrete sequence $0 = \lambda_0 < \lambda_1 \leq \lambda_2 \leq \ldots \to \infty$ of the eigenvalues of $\Delta_g$, counted with the corresponding multiplicities, is called the *spectrum* of $(M, g)$; we will denote it by $\mathrm{spec}(M, g)$ or $\mathrm{spec}(\Delta_g)$. If $\mathcal{H} \subseteq L^2(M, g)$ is a subspace of functions which is invariant under $\Delta_g$, we will denote the corresponding spectrum of eigenvalues by $\mathrm{spec}(\mathcal{H})$. Two compact Riemannian manifolds are said to be *isospectral* if their spectra coincide.

In the following proposition we formulate some version of a principle which was observed by C. Gordon [Go2] and is the key of the isospectrality proofs also for the manifolds constructed in [GW3] and [GGSWW].

Let $T$ be a compact abelian Lie group, i.e., a torus, and let $L : T \times M \to M$ be a smooth free action of $T$ on a compact Riemannian manifold $(M, g)$ by isometries. There is a unique Riemannian metric, denoted $g^T$, on the quotient manifold $T\backslash M$ such that the canonical projection $\pi_T : (M,g) \to (T\backslash M, g^T)$ becomes a Riemannian submersion.

1.2 PROPOSITION. *Let $(M,g)$ and $(M',g')$ be two closed Riemannian manifolds, and let $L : T \times M \to M$ and $L' : T \times M' \to M'$ be smooth free actions of a torus $T$ by isometries. Assume that for every closed subgroup $W$ of $T$ which is either $T$ itself or a subtorus of codimension $1$ in $T$, the following holds:*

(i) *The fibers of the projections $\pi_W : M \to W\backslash M$ and $\pi'_W : M' \to W\backslash M'$ are totally geodesic submanifolds of $(M,g)$ and $(M',g')$, respectively.*

(ii) *The Riemannian manifolds $(W\backslash M, g^W)$ and $(W\backslash M', g'^W)$ are isospectral.*

*Then $(M,g)$ and $(M',g')$ are isospectral.*

*Proof.* Write $T = \mathbb{R}^r/\mathcal{L}$, where $\mathcal{L}$ is a lattice of full rank in $\mathbb{R}^r$. For every $\bar{z}$ in $T$ choose $z \in \mathbb{R}^r$ such that $\bar{z} = z + \mathcal{L}$. Consider the unitary representation $\rho$ of $T$ on $L^2(M,g)$ defined by $\rho_{\bar{z}} f = f \circ L_{\bar{z}}$ for all $\bar{z} \in T$, $f \in L^2(M,g)$. Since $T$ is abelian, $L^2(M,g)$ decomposes as an orthogonal sum $\bigoplus_{\mu \in \mathcal{L}^*} \mathcal{H}_\mu$, where $\mathcal{L}^*$ is the dual lattice of $\mathcal{L}$, and

$$\mathcal{H}_\mu = \{f \in L^2(M,g) \mid \rho_{\bar{z}} f = e^{2\pi i \mu(z)} f \quad \text{for all } \bar{z} \in T\}.$$

The spaces $\mathcal{H}_\mu$ are invariant under $\Delta := \Delta_g$ because $T$ acts by isometries. Note that $\pi_T^*$ is a linear bijection from $L^2(T\backslash M, g^T)$ to $\mathcal{H}_0$. Since $\pi_T$ is a Riemannian submersion with totally geodesic fibers by (i), $\pi_T^*$ intertwines the corresponding Laplacians. Thus $\text{spec}(\mathcal{H}_0) = \text{spec}(T\backslash M, g^T)$. Assumption (ii) for $W = T$ implies, with the obvious analogous notations for $(M',g')$, that $\text{spec}(\mathcal{H}_0) = \text{spec}(\mathcal{H}'_0)$.

Now let $W$ be a closed subgroup of codimension $1$ in $T$. Denote by $\widetilde{W}$ the corresponding subspace of $\mathbb{R}^r$. Then $\pi_W^*$ is a linear bijection from $L^2(W\backslash M, g^W)$ to $\bigoplus_{\substack{\mu \in \mathcal{L}^* \\ \mu|\widetilde{W}=0}} \mathcal{H}_\mu = \mathcal{H}_0 \oplus \bigoplus_{\substack{\mu \in \mathcal{L}^* \\ \ker \mu = \widetilde{W}}} \mathcal{H}_\mu$. From assumptions (i) and (ii) for $W$ and the fact that $\text{spec}(\mathcal{H}_0) = \text{spec}(\mathcal{H}'_0)$ we conclude, using the same argument as before, that

$$\text{spec}\Big(\bigoplus_{\substack{\mu \in \mathcal{L}^* \\ \ker \mu = \widetilde{W}}} \mathcal{H}_\mu\Big) = \text{spec}\Big(\bigoplus_{\substack{\mu \in \mathcal{L}^* \\ \ker \mu = \widetilde{W}}} \mathcal{H}'_\mu\Big).$$

Since $\mathcal{L}^* \setminus \{0\}$ is the disjoint union of the sets $\{\mu \in \mathcal{L}^* \mid \ker \mu = \widetilde{W}\}$, where $W$ runs through the set of all subtori of codimension $1$ in $T$, we conclude $\text{spec}(M,g) = \text{spec}(M',g')$. □

1.3 *Definitions and Notation.* (i) For positive integers $m$ and $r$ let $\mathfrak{v} := \mathbb{R}^m$ and $\mathfrak{z} := \mathbb{R}^r$. Throughout this paper, we will assume that $\mathfrak{v}$ and $\mathfrak{z}$ are equipped with the euclidean standard metric. Moreover, $j$ will always denote a linear map from $\mathfrak{z}$ to $so(\mathfrak{v})$. With every such map $j$ we associate a two-step nilpotent Lie algebra $\mathfrak{g}_j = \mathfrak{v} \oplus \mathfrak{z}$ whose Lie bracket is defined by $\langle z, [x,y] \rangle = \langle j_z x, y \rangle$ for all $x, y \in \mathfrak{v}$, $z \in \mathfrak{z}$. Every two-step nilpotent Lie algebra arises in this way. Let $G_j$ be the simply connected Lie group with Lie algebra $\mathfrak{g}_j$. The Lie group exponential map $\exp : \mathfrak{g}_j \to G_j$ is a diffeomorphsim. We write the elements of $G_j$ as $(x, z) := \exp(x + z)$ with $x \in \mathfrak{v}$ and $z \in \mathfrak{z}$. With this notation, multiplication in $G_j$ is given by

$$(1) \qquad (x, z) \cdot (y, w) = (x + y, z + w + \tfrac{1}{2}[x,y]).$$

(ii) Let $S$ be a simply connected compact Lie group whose maximal torus is of dimension not less than $r$. Denote the Lie algebra of $S$ by $\mathfrak{s}$. Let $H \subset S$ be an $r$-dimensional closed abelian subgroup with Lie algebra $\mathfrak{h} \subset \mathfrak{s}$. Let $S$ be equipped with some fixed biinvariant metric $k$, and consider the corresponding scalar product on $\mathfrak{s}$. Choose a linear bijection $f : \mathfrak{z} \to \mathfrak{h}$ which is an isometry with respect to the induced metric on $\mathfrak{h}$ and the standard metric on $\mathfrak{z}$. Define

$$E := \mathfrak{v} \times S,$$

and define actions of $G_j$ on $E$ from the left and of $S$ on $E$ from the right by

$$(2) \qquad L_{(x,z)}(v, s) := \bigl(x + v, \exp(f(z + \tfrac{1}{2}[x,v])) \cdot s\bigr), \quad R_t(v, s) := (v, st).$$

Here $\exp$ denotes the Lie group exponential map from $\mathfrak{s}$ to $S$. Using the multiplication law (1), one easily verifies that $L$ is indeed a group action. Note that these actions of $G_j$ and $S$ commute. Let $D_j := L_{G_j} \cdot R_S$. Then $D_j$ is a transitive group of diffeomorphisms of $E$.

1.4 LEMMA. *If $\varphi \in D_j$ fixes $(0, e) \in \mathfrak{v} \times S = E$, then $\varphi_{*(0,e)}$ is of the form $(\mathrm{Id}, \psi)$, where $\psi = \mathrm{Ad}_h$ for some $h \in H$.*

*Proof.* Let $\varphi = L_{(x,z)} \cdot R_s$. If $\varphi$ fixes $(0, e)$ then $x = 0$ and $\exp(f(z)) = s^{-1}$. Thus $s \in H$, and $\varphi$ acts on $E$ by $(x, t) \mapsto (x, s^{-1}ts)$. Therefore $\varphi_{*(0,e)} = (\mathrm{Id}, \mathrm{Ad}_s^{-1})$. □

1.5 *Definition.* Define a metric $g$ on $T_{(0,e)}E = \mathfrak{v} \oplus \mathfrak{s}$ by letting $\mathfrak{v}$ and $\mathfrak{s}$ be orthogonal and requiring that the induced metric on $\mathfrak{v}$, resp. $\mathfrak{s}$, be the standard metric on $\mathfrak{v}$, resp. the scalar product corresponding to the biinvariant metric $k$. Let $g_j$ be the Riemannian metric on $E$ which is invariant under $D_j$ and coincides with $g$ on $T_{(0,e)}E$.

*Remark.* Since $k$ corresponds to a biinvariant metric on $S$, inner automorphisms of $\mathfrak{s}$ are orthogonal with respect to $k$. Thus Lemma 1.4 implies that $g_j$ is well-defined.

1.6 LEMMA. *For all $x \in \mathfrak{v}$ and $s \in S$, the submanifolds $\{x\} \times S$ and $\{x\} \times Hs$ are totally geodesic in $(E, g_j)$. The same holds for $\{x\} \times Ws$, where $W$ is any closed subgroup of $H$.*

*Proof.* Consider the involution $\sigma : E \to E$ given by $\sigma(x,s) = (-x,s)$. Obviously $\sigma_{*(0,e)}$ is an isometry of the metric $g$ on $T_{(0,e)}E$. Moreover, $\sigma \circ R_s = R_s \circ \sigma$ for all $s \in S$, and one easily checks by (2) that $\sigma \circ L_{(x,z)} = L_{(-x,z)} \circ \sigma$ for all $(x,z) \in G_j$. These equations imply that $\sigma^* g_j$ is invariant under $D_j$ and therefore equal to $g_j$. Thus $\sigma$ is an isometry of $(E, g_j)$; hence its fixed point set $\{0\} \times S$ is totally geodesic. Note that the metric induced on $\{0\} \times S$ by $g_j$ corresponds to the metric $k$ on $S$ under the canonical bijection. Since $H$ is totally geodesic in $S$ with respect to this metric, the set $\{0\} \times H$ is totally geodesic in $(E, g_j)$. Moreover, $\{0\} \times H$ is a flat torus. Hence any subtorus $\{0\} \times W$ corresponding to a closed subgroup $W$ of $H$ is totally geodesic in $\{0\} \times H$ and thus also in $(E, g_j)$. For $x \in \mathfrak{v}$ and $s \in S$, the isometry $L_{(x,0)} \circ R_s$ carries $\{0\} \times S$, $\{0\} \times H$, and $\{0\} \times W$ to $\{x\} \times S$, $\{x\} \times Hs$, and $\{x\} \times Ws$, respectively, whence the statement follows. □

The following lemma is immediate from (2). Part (i) describes the $g_j$-orthogonal splitting of the tangent space $T_{(x,s)}E$ at an arbitrary point $(x,s) \in E$ into horizontal and vertical parts with respect to the projection $E = \mathfrak{v} \times S \to \mathfrak{v}$.

1.7 LEMMA. (i) *For every $(x,s) \in E$, the $g_j$-isometry $(L_{(x,0)} \circ R_s)_* : T_{(0,e)}E = \mathfrak{v} \oplus \mathfrak{s} \to T_{(x,s)}E$ maps $v \in \mathfrak{v}$ to $v^j_{(x,s)} := (v, r_{s*}f(\frac{1}{2}[x,v]))$ and $u \in \mathfrak{s}$ to $u^j_{(x,s)} := (0, r_{s*}u)$, where $r_s$ denotes right multiplication by $s$ in $S$.*

(ii) *For $v \in \mathfrak{v}$ and $h \in \mathfrak{h} \subset \mathfrak{s}$, the vector fields $v^j$ and $h^j$ are invariant under the group $D_j = L_{G_j} \cdot R_S$. For a vector $u$ in the orthogonal complement $\mathfrak{u}$ of $\mathfrak{h}$ in $\mathfrak{s}$, the vector field $u^j$ is invariant under $R_S$, but in general not under $L_{G_j}$. However, the distribution $\mathfrak{u}_{(x,s)} = \{u^j_{(x,s)} \mid u \in \mathfrak{u}\}$ is invariant under both actions.*

Note that for $v \in \mathfrak{v}$, the vector field $v^j$ on $E$ depends on $j$ since $[\,,\,]$ does so, while $u^j$ is actually independent of $j$ for $u \in \mathfrak{s}$. We use the superscript just to distinguish between the vector field $u^j$ on $E$ and the vector $u \in \mathfrak{s}$.

1.8 *Definition.* Let $C := S^1(\mathfrak{v}) \times S \subset E$, where $S^1(\mathfrak{v})$ is the unit sphere in $\mathfrak{v}$. By $g_j$ we will in the following denote the above metric on $E$ as well as its restriction to $C$.

1.9 *Remark.* Concerning the geometry of the manifolds $(C, g_j)$, note that the canonical projection from $(C, g_j)$ to its first factor $S^1(\mathfrak{v})$, endowed with the round sphere metric, is a Riemannian submersion with totally geodesic fibers which are isometric to the compact Lie group $(S, k)$. However, the inclusion $S^1(\mathfrak{v}) \to S^1(\mathfrak{v}) \times \{e\} \subset (C, g_j)$ is not isometric, nor is its image orthogonal to the fibers $\{x\} \times S$ of the above submersion.

## 2. The isospectral metrics

2.1 *Definition* ([GW3]). Two linear maps $j, j' : \mathfrak{z} \to \mathrm{so}(\mathfrak{v})$ are called *isospectral*, denoted $j \sim j'$, if for every $z \in \mathfrak{z}$ there exists an orthogonal map $A_z : \mathfrak{v} \to \mathfrak{v}$ such that $j'_z = A_z j_z A_z^{-1}$.

2.2 THEOREM. *If $j \sim j'$ then $(C, g_j)$ and $(C, g_{j'})$ are isospectral.*

*Proof.* Let $\mathcal{L}$ be the set of those $z \in \mathfrak{z}$ for which $L_{(0,z)} = \mathrm{Id}$. Thus $\mathcal{L} = f^{-1}\{h \in \mathfrak{h} \mid \exp h = e\}$ in the notation 1.3(ii). Since $H$ is a torus, $\mathcal{L}$ is a lattice of full rank in $\mathfrak{z}$. The linear isometry $f : \mathfrak{z} \to \mathfrak{h}$ induces an isometry $\bar{f}$ between the flat tori $T := \mathfrak{z}/\mathcal{L}$ and $H$. The torus $T$ acts on $E$ by $L_{\bar{z}} := L_{(0,z)}$, where $\bar{z} = z + \mathcal{L} \in T$. Obviously this action is smooth, free, and by isometries with respect to $g_j$. The same holds for its restriction to the submanifold $C$ which is invariant under the action of $T$. The orbits are of the form $\{v\} \times Hs \subset C$ with $v \in S^1(\mathfrak{v})$ and $s \in S$, and the orbits of the action of any closed subgroup $W$ of $T$ have the form $\{v\} \times Us$, where $U = \bar{f}(W)$ is the corresponding subtorus of $H$. By Lemma 1.6 these orbits are totally geodesic in $(E, g_j)$ and thus also in $(C, g_j)$.

Since these facts are true also with respect to $g_{j'}$, we are in the situation of Proposition 1.2 with condition (i) satisfied. Our goal now is to show that $j \sim j'$ implies that also condition (ii) of the proposition is satisfied. We will actually show that the various pairs of quotients $(W \backslash C, g_j^W)$ and $(W \backslash C, g_{j'}^W)$ are not only isospectral here, but even isometric.

First consider $W = T$. Denote the Lie bracket of $\mathfrak{g}_{j'}$ by $[\,,\,]'$. Then for all $x, v \in \mathfrak{v}$ and $s \in S$ the difference between the vectors $v^j_{(x,s)}$ and $v^{j'}_{(x,s)} \in T_{(x,s)}E$ equals $(0, r_{s*}f(\frac{1}{2}[x,v] - \frac{1}{2}[x,v]'))$ and thus is tangent to the $L_T$-orbit through $(x, s)$. Lemma 1.7 now implies that $\mathrm{Id} : E \to E$ induces an isometry between the quotient manifolds $(T \backslash E, g_j^T)$ and $(T \backslash E, g_{j'}^T)$. By restriction, $\mathrm{Id} : C \to C$ induces the required isometry between $(T \backslash C, g_j^T)$ and $(T \backslash C, g_{j'}^T)$.

Now let $W \subset T$ be a closed subgroup of codimension 1. Denote by $\widetilde{W}$ the corresponding subspace of $\mathfrak{z}$, and choose $z \in \mathfrak{z} \setminus \{0\}$ such that $z \perp \widetilde{W}$. Since $j \sim j'$, there exists an orthogonal map $A_z : \mathfrak{v} \to \mathfrak{v}$ such that $j'_z = A_z j_z A_z^{-1}$. Define a diffeomorphism $F_z : E \to E$ by $F_z(x, s) = (A_z x, s)$. For $u \in \mathfrak{s}$, the vector field $u^j = u^{j'}$ from Lemma 1.7 is obviously carried to itself by $F_{z*}$. Since the tangent space to the $L_W$-orbit through $(x, s) \in E$ is just $\{u^j_{(x,s)} \mid u \in f(\widetilde{W})\}$, the map $F_z$ carries $L_W$-orbits to $L_W$-orbits. Moreover, for every $v \in \mathfrak{v}$ we have

$$F_{z*}|_{(x,s)}(v^j) = \left(A_z v, r_{s*}f(\tfrac{1}{2}[x,v])\right) \quad \text{and}$$

$$(A_z v)^{j'}|_{F_z(x,s)} = \left(A_z v, r_{s*}f(\tfrac{1}{2}[A_z x, A_z v]')\right).$$

The choice of $A_z$ implies that $\langle z, [x,v] - [A_zx, A_zv]'\rangle = \langle j_zx, v\rangle - \langle j'_z A_zx, A_zv\rangle = 0$. Thus the difference between $F_{z*}|_{(x,s)}(v^j)$ and $(A_zv)^{j'}|_{F_z(x,s)}$ is orthogonal to $(0, r_{s*}f(z))$ and therefore tangent to the $L_W$-orbit through $F_z(x,s)$. Using Lemma 1.7 we conclude that $F_z$ induces an isometry between the quotients $(W\backslash E, g_j^W)$ and $(W\backslash E, g_{j'}^W)$. Since $F_z$ preserves $C$, this isometry restricts to the required isometry between $(W\backslash C, g_j^W)$ and $(W\backslash C, g_{j'}^W)$. □

2.3 *Remark.* In our notation, the bounded isospectral manifolds constructed in [GW3] are of the form $(B^1(\mathfrak{v}) \times H, g_j)$, resp. $(B^1(\mathfrak{v}) \times H, g_{j'})$, where $B^1(\mathfrak{v})$ denotes the unit ball in $\mathfrak{v}$. Similarly, the isospectral manifolds in [GGSWW] are of the form $(S^1(\mathfrak{v}) \times H, g_j)$, resp. $(S^1(\mathfrak{v}) \times H, g_{j'})$.

2.4 *Remarks.* (i) If there exists an orthogonal map $A: \mathfrak{v} \to \mathfrak{v}$ such that $j'_z = Aj_zA^{-1}$ for all $z \in \mathfrak{z}$ then the isospectrality condition from Definition 2.1 is trivially satisfied. In this situation the manifolds $(C, g_j)$ and $(C, g_{j'})$ are isometric; an isometry from the first to the latter is given by $(x,s) \mapsto (Ax, s)$. For example, this is always the case if $\dim \mathfrak{z} = 1$.

(ii) However, we will show that in case $\dim \mathfrak{z} = 2$ and $\dim \mathfrak{v} \geq 5$, any generic linear map $j: \mathfrak{z} \to so(\mathfrak{v})$ is contained in a continuous isospectral family of maps $j(t)$ such that the manifolds $(C, g_{j(t)})$ are not pairwise isometric. Here, *generic* means that $j$ is an element of a certain Zariski open subset of the space of all linear maps from $\mathfrak{z}$ to $so(\mathfrak{v})$. Note that a similar genericity result, although excluding the case $\dim \mathfrak{v} = 6$, was obtained in [GW3, Theorem 2.2 and Proposition 1.4] and [GGSWW, Propositions 9 and 10] with respect to the isospectral manifolds considered there. That result was stronger than the one we prove here in the sense that it included the existence of multi-parameter families of isospectral metrics. On the other hand, our result includes the case $\dim \mathfrak{v} = 6$ in which the previously used arguments had failed.

(iii) For $\dim \mathfrak{v} \leq 4$ and $\dim \mathfrak{z} = 2$, it is not hard to show that, as mentioned in [GW3], the condition $j \sim j'$ can be satisfied only trivially (in the sense of (i) above). For $\dim \mathfrak{v} \leq 4$ and $\dim \mathfrak{z} \geq 3$ it is still possible to show that there are no *continuous* families of nontrivially isospectral maps $j(t): \mathfrak{z} \to so(\mathfrak{v})$ (although a *pair* of such maps is given in [GW3, Example 1.10]). Thus $\dim \mathfrak{v} = 5$ and $\dim \mathfrak{z} = 2$ is indeed the minimal choice of dimensions for which nontrivial continuous isospectral families $j(t): \mathfrak{z} \to so(\mathfrak{v})$ exist.

(iv) The lowest dimensional compact, simply connected Lie group $S$ whose maximal torus has dimension $r$ is the direct product of $r$ copies of SU(2). For this choice of $S$, the manifold $C$ is diffeomorphic to the product of $S^{m-1}$ with $r$ copies of $S^3$. By the previous remark, the lowest dimensional example of this form on which our method produces nontrivial, continuous families of isospectral metrics is $S^4 \times S^3 \times S^3$.

2.5 *Example.* This example is a modified version of Example 2.3 from [GW3]. The modification mainly consists in lowering the dimension of $\mathfrak{v}$ by one. Consider $\mathfrak{v} = \mathbb{R}^5$ and $\mathfrak{z} = \mathbb{R}^2$. Let $a, b \in \mathrm{so}(\mathfrak{v})$ be given as matrices with respect to the standard basis by

$$a = \begin{pmatrix} 0 & -a_1 & 0 & 0 & 0 \\ a_1 & 0 & 0 & 0 & 0 \\ 0 & 0 & 0 & -a_2 & 0 \\ 0 & 0 & a_2 & 0 & 0 \\ 0 & 0 & 0 & 0 & 0 \end{pmatrix}, \quad b = \begin{pmatrix} 0 & 0 & b_{12} & 0 & b_{13} \\ 0 & 0 & 0 & 0 & 0 \\ -b_{12} & 0 & 0 & 0 & b_{23} \\ 0 & 0 & 0 & 0 & 0 \\ -b_{13} & 0 & -b_{23} & 0 & 0 \end{pmatrix}.$$

Define a linear map $j_{a,b} : \mathfrak{z} \to \mathrm{so}(\mathfrak{v})$ by letting $j_{a,b}(z_1) = a$ and $j_{a,b}(z_2) = b$, where $\{z_1, z_2\}$ is the standard basis of $\mathfrak{z}$. The characteristic polynomial $\det(\lambda \mathrm{Id} - (sa + ub))$ is equal to $\lambda^5 + (s^2(a_1^2 + a_2^2) + u^2(b_{12}^2 + b_{13}^2 + b_{23}^2))\lambda^3 + (s^4 a_1^2 a_2^2 + s^2 u^2(a_1^2 b_{23}^2 + a_2^2 b_{13}^2))\lambda$. Thus $j_{a,b} \sim j_{a,b'}$ if and only if

$$b_{12}^2 + b_{13}^2 + b_{23}^2 = b_{12}'^2 + b_{13}'^2 + b_{23}'^2 \quad \text{and} \quad a_1^2 b_{23}^2 + a_2^2 b_{13}^2 = a_1^2 b_{23}'^2 + a_2^2 b_{13}'^2.$$

Now assume $0 < a_1 < a_2$ and $b_{12}, b_{13}, b_{23} \geq 0$. Then with

$$b_{12}(t) := \sqrt{b_{12}^2 + t(a_2^2 - a_1^2)}, \quad b_{13}(t) := \sqrt{b_{13}^2 + ta_1^2}, \quad b_{23}(t) := \sqrt{b_{23}^2 - ta_2^2},$$

the maps $j_{a,b}$ and $j_{a,b(t)}$ satisfy conditions (i) and (ii) for every

$$t \in \left[\max\left\{\frac{-b_{12}^2}{a_2^2 - a_1^2}, \frac{-b_{13}^2}{a_1^2}\right\}, \frac{b_{23}^2}{a_2^2}\right].$$

This interval has nonzero length if $b_{23} > 0$ or $b_{12}, b_{13} > 0$, and $j(t) := j_{a,b(t)}$ is a continuous family of pairwise isospectral linear maps. This family is nontrivial; i.e., the maps $j(t)$ are not pairwise equivalent in the sense of Remark 2.4(i). This can be seen, for example, by noting that the maps $(j_{z_1}^2 + j_{z_2}^2)(t) = a^2 + b(t)^2$ are not pairwise similar, since routine computation shows that their determinant is a linear polynomial in $t$ with nonzero leading coefficient $a_1^4 a_2^4(a_1^2 - a_2^2)$. In the next section, we will interpret the nonsimilarity of the various $(j_{z_1}^2 + j_{z_2}^2)(t)$ geometrically in terms of the scalar curvature of the isospectral manifolds $(C, g_{j(t)})$ (see formula (10) and Proposition 3.5).

## 3. Scalar curvature and a sufficient condition for nonisometry

In this section we will establish a sufficient condition for two manifolds $(C, g_j)$ and $(C, g_{j'})$ (with $j \sim j'$) to be nonisometric. We will show that the difference in geometry can in general be detected by studying the behaviour of the scalar curvature. The manifolds $(C, g_j)$ are in general not locally homogeneous, and the scalar curvature function is in general nonconstant on a given manifold $(C, g_j)$. In Proposition 3.5, we will formulate a condition on two isospectral linear maps $j, j' : \mathfrak{z} \to \mathrm{so}(\mathfrak{v})$ which will imply that $(C, g_j)$ and $(C, g_{j'})$ have different sets of critical values of the scalar curvature function. In Section four we will then show that in case $\dim \mathfrak{z} = 2$ and $\dim \mathfrak{v} \geq 5$, any generic linear map $j : \mathfrak{z} \to \mathrm{so}(\mathfrak{v})$ admits a continuous isospectral deformation $j(t)$ such that the maps $j(t)$ pairwise satisfy this condition.

We need some preparations in order to compute the scalar curvature of $(C, g_j)$. In the following, let $G_j$ be endowed with the left invariant metric $m_j$ induced by the standard scalar product on $\mathfrak{g}_j = \mathfrak{v} \oplus \mathfrak{z}$.

3.1 LEMMA. *The map $\ell : G_j \ni (x, z) \mapsto L_{(x,z)}(0, e) = (x, e^{f(z)}) \in L_{G_j}(0, e) \subset (E, g_j)$ is a local isometry. Moreover, the submanifold $L_{G_j}(0, e)$ is totally geodesic in $(E, g_j)$; consequently, this holds also for its right translates under $R_S$.*

*Proof.* Two elements of $G_j$ have the same image under $\ell$ if and only if they differ by an element of $\exp \mathcal{L}$, where $\mathcal{L}$ is the central lattice $f^{-1}\{h \in \mathfrak{h} \mid \exp h = e\} \subset \mathfrak{z}$ as in the proof of Theorem 2.2. Thus $\ell$ is a local diffeomorphism. For $v \in \mathfrak{v}$, the left invariant vector field on $G_j$ induced by $v$ is mapped by $\ell_*$ to $v^j|_{L_{G_j}(0,e)}$, and for $z \in \mathfrak{z}$, the corresponding left invariant vector field on $G_j$ is mapped by $\ell_*$ to $h^j|_{L_{G_j}(0,e)}$, where $h = f(z)$. Lemma 1.7 now implies that $\ell$ is a local isometry from $(G_j, m_j)$ to $L_{G_j}(0, e)$.

To see that $L_{G_j}(0, e)$ is totally geodesic, first note that its tangent space at $(0, e)$ is $\mathfrak{v} \oplus \mathfrak{h} \subset \mathfrak{v} \oplus \mathfrak{s} = T_{(0,e)}E$, and $L_{G_j}(0, e)$ is the integral manifold through $(0, e)$ of the distribution $\{v^j + h^j \mid v \in \mathfrak{v}, h \in \mathfrak{h}\}$. In particular, for all $x, v \in \mathfrak{v}$ and every $u \in \mathfrak{u}$, where $\mathfrak{u}$ denotes the orthogonal complement of $\mathfrak{h}$ in $\mathfrak{s}$, we have

$$(3) \qquad \langle [x^j, v^j]_{(0,e)}, u \rangle = 0.$$

Extend every $u \in \mathfrak{u} \subset T_{(0,e)}E$ to a vector field $\tilde{u}$ on $E$ whose flow $F_{\tilde{u}}^t$ consists of the right translations $R_{e^{tu}}$ from 1.3(ii), and consider the vector field $v^j$ associated with $v \in \mathfrak{v}$ as in Lemma 1.7. One easily checks that the flow $F_{v^j}^t$ commutes with the maps $R_{e^{tu}}$; thus

$$(4) \qquad [v^j, \tilde{u}] = 0.$$

Moreover, for $x, v \in \mathfrak{v}$ the function $\langle x^j, v^j \rangle$ is constant on $E$ by Lemma 1.7, and $\langle x^j, \tilde{u} \rangle = \langle v^j, \tilde{u} \rangle = 0$ on $E$. By the Koszul formula for the Levi-Cività connection $\nabla$ on $(E, g_j)$, this implies together with (3) and (4) that

$$(5) \qquad \langle \nabla_x v^j, u \rangle = 0.$$

We know that $\nabla_h h^j = 0$ for $h \in \mathfrak{h}$, since $h^j$ is tangent to the flat submanifold $\{0\} \times H$ which is totally geodesic by Lemma 1.6. Finally, for all $v \in \mathfrak{v}, h \in \mathfrak{h}$, and $u \in \mathfrak{u}$, we have $\langle \nabla_v h^j, u \rangle = \langle \nabla_{\sigma_* v} \sigma_* h^j, \sigma_* u \rangle = \langle \nabla_{-v} h^j, u \rangle$, where $\sigma : E \ni$ $(x,s) \mapsto (-x,s) \in E$ is the involutive isometry from the proof of 1.6. Thus $\langle \nabla_v h^j, u \rangle = 0$, and similarly $\langle \nabla_h v^j, u \rangle = 0$. In summary, $\nabla_y y^j \in \mathfrak{v} \oplus \mathfrak{h}$ for every $y \in \mathfrak{v} \oplus \mathfrak{h}$, whence $L_{G_j}(0,e)$ is totally geodesic. □

Denote the scalar curvatures of $(E, g_j)$, $(G_j, m_j)$, and $(S, k)$ by $\mathrm{scal}^{E,j}$, $\mathrm{scal}^{G_j}$, and $\mathrm{scal}^S$, respectively. Since all three manifolds are homogeneous, their scalar curvatures are constant.

3.2 LEMMA.
$$\mathrm{scal}^{E,j} = \mathrm{scal}^{G_j} + \mathrm{scal}^S.$$

*Proof.* We compute the scalar curvature in the point $(0,e) \in E$. Recall the orthogonal splitting $T_{(0,e)}E = \mathfrak{v} \oplus \mathfrak{h} \oplus \mathfrak{u}$ from the proof 3.1. Since $\mathfrak{v} \oplus \mathfrak{h}$ and $\mathfrak{h} \oplus \mathfrak{u} = \mathfrak{s}$ are the tangent spaces to $L_{G_j}(0,e)$ and $\{0\} \times S$ which are both totally geodesic by 3.1 and 1.6, and since the sectional curvature of planes tangent to $\mathfrak{h}$ is zero, it suffices to show that

(6) $$K(v, u) = 0$$

for all $v \in \mathfrak{v}$ and $u \in \mathfrak{u}$, where $K$ denotes the sectional curvature on $(E, g_j)$.

Using again the involutive isometry $\sigma : E \ni (x,s) \mapsto (-x,s) \in E$, we have $\langle \nabla_u v^j, w \rangle = \langle \nabla_u(-v^j), w \rangle = 0$ for every $w \in \mathfrak{s}$. Thus $\nabla_u v^j \in \mathfrak{v}$. But the same argument which led to (5) also shows $\langle \nabla_u v^j, x \rangle = 0$ for every $x \in \mathfrak{v}$. Hence $\nabla_u v^j = 0$, and by Lemma 1.7,

(7) $$\nabla_{u^j} v^j = 0 \quad \text{for all } u \in \mathfrak{u}, v \in \mathfrak{v}.$$

Since $u^j$ equals $\tilde{u}$ along the curve $t \mapsto (tv, e)$, where $\tilde{u}$ is the vector field from the proof of 3.1, we get $\nabla_v \nabla_{\tilde{u}} v^j = 0$. By (4) we also know $\nabla_{[\tilde{u}, v^j]} v^j = 0$ and thus $R(u,v)v = \nabla_u \nabla_{v^j} v^j$. Note that $\nabla_{v^j} v^j = 0$ by Lemma 3.1 and the fact that $\nabla_v v = -{}^t\mathrm{ad}_v v = 0$ for the corresponding left invariant vector field on $(G_j, m_j)$, since $v \perp [\mathfrak{g}_j, \mathfrak{g}_j]$. Hence $R(u,v)v = 0$, and the statement follows. □

We can now compute the scalar curvature of the submanifold $(C, g_j)$ of $(E, g_j)$:

3.3 LEMMA. *The scalar curvature $\mathrm{scal}^{C,j}$ of $(C, g_j)$ at the point $(v,s) \in C$ is given by*

$$\mathrm{scal}^{C,j}_{(v,s)} = \mathrm{scal}^{G_j} + \mathrm{scal}^S + (m-1)(m-2) - \mathrm{Ric}^{G_j}(v,v),$$

*where $m = \dim(\mathfrak{v})$, and $\mathrm{scal}^{G_j}$ and $\mathrm{Ric}^{G_j}$ denote the scalar curvature and the Ricci tensor of $(G_j, m_j)$.*

*Proof.* Let $N$ be the outward unit normal field to $(C, g_j)$ in $(E, g_j)$. A standard computation involving only the Gauss equations and the fact that

$(C, g_j)$ is a codimension one submanifold of $(E, g_j)$ shows that

$$(8) \quad \mathrm{scal}_p^{C,j} = \mathrm{scal}^{E,j} - 2\mathrm{Ric}(N_p, N_p) + \left(\mathrm{tr}(\nabla N_{|T_pC})\right)^2 - \|\nabla N_{|T_pC}\|^2$$

for all $p \in C$, where $\nabla$ and Ric denote the Levi-Cività connection and the Ricci tensor of $(E, g_j)$, and $\|.\|$ denotes the euclidean norm on tensors. Note that for $p = (v, s)$, the tangent space $T_pC$ is the orthogonal sum of $\{x^j_{(v,s)} \mid x \in \mathfrak{v}, x \perp v\}$ and $\{u^j_{(v,s)} \mid u \in \mathfrak{s}\}$. Thus $N_{(v,s)} = v^j_{(v,s)}$. Moreover, one easily checks that for each $x \perp v$ in $\mathfrak{v}$ and each $u \in \mathfrak{s}$, we have

$$\nabla_{x^j_{(v,s)}} N = x^j_{(v,s)} + \nabla_{x^j_{(v,s)}} v^j, \quad \nabla_{u^j_{(v,s)}} N = \nabla_{u^j_{(v,s)}} v^j.$$

Recall that $L_{G_j}$-orbits are totally geodesic in $(E, g_j)$ by Lemma 3.1. Since $x^j$, $v^j$, and $h^j$ for $h \in \mathfrak{h}$ are tangent to these, we get for all $x \in \mathfrak{v}$ and $h \in \mathfrak{h}$, with $z := f^{-1}(h) \in \mathfrak{z}$,

$$\nabla_{x^j_{(v,s)}} N = x^j_{(v,s)} + (\nabla_x^{G_j} v)^j_{(v,s)}, \quad \nabla_{h^j_{(v,s)}} N = (\nabla_z^{G_j} v)^j_{(v,s)},$$

whereas $\nabla_{u^j} N = 0$ for $u \in \mathfrak{u} = \mathfrak{h}^\perp \cap \mathfrak{s}$ by (7). We have (see e.g. [Eb])

$$(9) \quad \nabla_x^{G_j} v = \frac{1}{2}[x, v], \quad \nabla_z^{G_j} v = -\frac{1}{2} j_z v$$

for all $x, v \in \mathfrak{v}$ and $z \in \mathfrak{z}$. Thus $\nabla^{G_j} v$ sends $\mathfrak{v}$ to $\mathfrak{z}$ and vice versa; hence

$$\mathrm{tr}(\nabla N_{|T_{(v,s)}C}) = m - 1, \quad \|\nabla N_{|T_{(v,s)}C}\|^2 = m - 1 + \|\nabla^{G_j} v\|^2.$$

Moreover,

$$(10) \quad \mathrm{Ric}^{G_j}|_\mathfrak{v} = \frac{1}{2} \sum_{i=1}^r j_{z_i}^2,$$

where $\{z_1, \ldots, z_r\}$ is an orthonormal basis of $\mathfrak{z}$ (see [Eb]). One routinely derives from (9) and (10) that

$$\|\nabla^{G_j} v\|^2 = -\mathrm{Ric}^{G_j}(v, v).$$

Finally note that

$$\mathrm{Ric}(v^j, v^j) = \mathrm{Ric}^{G_j}(v, v)$$

for all $v \in \mathfrak{v}$. This follows using again Lemma 3.1 and the fact that the sectional curvatures $K(v, u)$ vanish by (6) whenever $v \in \mathfrak{v}$ and $u \in \mathfrak{u} = \mathfrak{h}^\perp \cap \mathfrak{s}$. The statement is now immediate from (8) and Lemma 3.2. □

The scalar curvature function $\mathrm{scal}^{C,j}$ can thus be viewed as a quadratic function in the first coordinate $v$ of $(v, s) \in S^1(\mathfrak{v}) \times S = C$. Consequently, we have:

3.4 COROLLARY. *The set of critical points of the scalar curvature on $(C, g_j)$ is given by $\{(v, s) \in C \mid v$ is an eigenvector of $\operatorname{Ric}^{G_j}|_{\mathfrak{v}}\}$, and the set of critical values is $\{\operatorname{scal}^{G_j} + \operatorname{scal}^S + (m-1)(m-2) - \lambda \mid \lambda$ is an eigenvalue of $\operatorname{Ric}^{G_j}|_{\mathfrak{v}}\}$.*

We can now establish a sufficient condition for nonisometry of two metrics $g_j$ and $g_{j'}$ on $C$. Note that by the discussion in 2.5, this condition is satisfied for the families given there.

3.5 PROPOSITION. *Let $j \sim j'$. If $\operatorname{Ric}^{G_j}|_{\mathfrak{v}}$ and $\operatorname{Ric}^{G_{j'}}|_{\mathfrak{v}}$ have different sets of eigenvalues, then $\operatorname{scal}^{C,j}$ and $\operatorname{scal}^{C,j'}$ have different sets of critical values. In particular, $(C, g_j)$ and $(C, g_{j'})$ are not isometric. More generally, if $\operatorname{Ric}^{G_j}|_{\mathfrak{v}}$ and $\operatorname{Ric}^{G_{j'}}|_{\mathfrak{v}}$ do not have the same collection of eigenvalues, counted with multiplicities (i.e., if $\operatorname{Ric}^{G_j}|_{\mathfrak{v}}$ and $\operatorname{Ric}^{G_{j'}}|_{\mathfrak{v}}$ are not similar), then $(C, g_j)$ and $(C, g_{j'})$ are not isometric.*

*Proof.* The first statement follows from Corollary 3.4 by noting that $j \sim j'$ implies $\operatorname{scal}^{G_j} = \operatorname{scal}^{G_{j'}}$. In fact, by [Eb] we have $\operatorname{Ric}^{G_j}(z, w) = -\frac{1}{4}\operatorname{tr}(j_z j_w)$ for all $z, w \in \mathfrak{z}$. Together with (10) we get $\operatorname{scal}^{G_j} = \frac{1}{4}\sum_{i=1}^r \operatorname{tr}(j_{z_i}^2)$, where $\{z_1, \ldots, z_r\}$ is an orthonormal basis of $\mathfrak{z}$. But if $j \sim j'$ then obviously $\operatorname{tr}(j_z^2) = \operatorname{tr}(j_z'^2)$ for each $z \in \mathfrak{z}$, and thus $\operatorname{scal}^{G_j} = \operatorname{scal}^{G_{j'}}$ as claimed.

The last statement follows from the fact that the multiplicity $n_\lambda$ of an eigenvalue $\lambda$ of $\operatorname{Ric}^{G_j}|_{\mathfrak{v}}$ can be read off from the dimension $\dim S + (n_\lambda - 1)$ of the set of critical points in which $\operatorname{scal}^{C,j}$ assumes the critical value $\operatorname{scal}^{G_j} + \operatorname{scal}^S + (m-1)(m-2) - \lambda$. □

3.6 *Remark.* Since the *total* scalar curvature is a spectral invariant, $(C, g_j)$ and $(C, g_{j'})$ must have the same total scalar curvature if $j \sim j'$. This follows also from 3.3 by direct computation; note that $\int_{S^1(\mathfrak{v})} \operatorname{Ric}^{G_j}(v, v)\, dv$ equals $\operatorname{vol}(S^{m-1}) \cdot \frac{1}{m} \cdot \operatorname{tr}(\operatorname{Ric}^{G_j}|_{\mathfrak{v}})$ which coincides for $j \sim j'$ by formula (10).

## 4. Generic families

Throughout this section we assume $\dim \mathfrak{z} = 2$. As mentioned in Remark 2.4(iii), nontrivial pairs of isospectral maps $j : \mathbb{R}^2 \to \mathfrak{so}(m)$ do not exist if $m \leq 4$. However, we will show in this section that in case $m \geq 5$, any generic linear map $j : \mathbb{R}^2 \to \mathfrak{so}(m)$ is contained in a nontrivial, continuous isospectral family of linear maps $j(t) : \mathbb{R}^2 \to \mathfrak{so}(m)$ which pairwise satisfy the nonisometry condition from Proposition 3.5. The idea is to construct the families $j(t)$ as the (locally defined) flow lines of a smooth vector field on the space of all linear maps from $\mathbb{R}^2$ to $\mathfrak{so}(m)$.

Taking into account Theorem 2.2 and Proposition 3.5, the following proposition not only completes the proof of the Main Theorem stated in the introduction, but also shows that there is an abundant class of examples to it.

4.1 PROPOSITION. *Let $m \geq 5$, and denote by $\mathcal{J}$ the space of all linear maps from $\mathbb{R}^2$ to $\mathrm{so}(m)$. Then there exists a Zariski open subset $U \subset \mathcal{J}$ such that every $j \in U$ is contained in a continuous isospectral family $j(t)$ with $j(0) = j$ such that the function $t \mapsto \|\mathrm{Ric}^{G_{j(t)}}|_{\mathfrak{v}}\|^2$ has nonzero derivative at $t = 0$. In particular, for every $j \in U$ there exists $\varepsilon(j) > 0$ such that for $-\varepsilon(j) < t < \varepsilon(j)$, the isospectral manifolds $(C, g_{j(t)})$ have pairwise different sets of critical values of the scalar curvature function.*

*Proof.* We write $j \in \mathcal{J}$ in the form $j = (j_1, j_2)$, where $j_1 := j_{z_1}$, $j_2 := j_{z_2}$, and $\{z_1, z_2\}$ is an orthonormal basis of $\mathfrak{z} = \mathbb{R}^2$. The isospectrality condition $j \sim j'$ says that for all $s, t \in \mathbb{R}$, the skew-symmetric map $sj_1 + tj_2$ has the same collection of eigenvalues, counted with multiplicities, as $sj'_1 + tj'_2$. As mentioned in the proof of [GW3, Theorem 2.2], it is easy to see that this is equivalent to the condition $\mathrm{tr}((sj_1 + tj_2)^k) = \mathrm{tr}((sj'_1 + tj'_2)^k)$ for all even numbers $k$ in $\{1, \ldots, m\}$, or equivalently, for all $k \in \mathbb{N}$. Let $a, b \in \mathbb{N}_0$ such that $a + b > 0$. Then the coefficient at the monomial $s^a t^b$ in $\mathrm{tr}((sj_1 + tj_2)^{a+b})$ equals

$$\tag{11} p_{a,b}(j) := \sum_{\sigma \in \mathfrak{S}_{a,b}} \mathrm{tr}(j_{\sigma(1)} \cdots j_{\sigma(a+b)}),$$

where $\mathfrak{S}_{a,b}$ denotes the set of all maps $\sigma : \{1, \ldots, a+b\} \to \{1, 2\}$ which satisfy $\#\sigma^{-1}(1) = a$, $\#\sigma^{-1}(2) = b$. Thus we get

$$\tag{12} j \sim j' \iff p_{a,b}(j) = p_{a,b}(j') \text{ for all } a, b \in \mathbb{N}_0 \text{ with } a + b > 0.$$

On the other hand, we consider the polynomial

$$j \mapsto \|\mathrm{Ric}^{G_j}|_{\mathfrak{v}}\|^2 = \|\tfrac{1}{2}(j_1^2 + j_2^2)\|^2 = \tfrac{1}{4}\big(\mathrm{tr}(j_1^4) + 2\mathrm{tr}(j_1^2 j_2^2) + \mathrm{tr}(j_2^4)\big).$$

The condition $j \sim j'$ implies $\mathrm{tr}(j_1^4) = \mathrm{tr}(j'^4_1)$ and $\mathrm{tr}(j_2^4) = \mathrm{tr}(j'^4_2)$. Thus in case $j \sim j'$ we have

$$\tag{13} \|\mathrm{Ric}^{G_j}|_{\mathfrak{v}}\|^2 - \|\mathrm{Ric}^{G_{j'}}|_{\mathfrak{v}}\|^2 = \tfrac{1}{2}\big(q(j) - q(j')\big),$$

where

$$\tag{14} q(j) := \mathrm{tr}(j_1^2 j_2^2).$$

By Lemma 4.2 below, the algebraic vector field $Y$ on $\mathcal{J}$, given by $Y(j) = (j_1^3 j_2 - j_2 j_1^3, 0)$, is everywhere orthogonal to the gradient of each of the polynomials $p_{a,b}$, while it is not everywhere orthogonal to the gradient of $q$. The first property implies by (12) that the flow lines $F_Y^t(j)$ consist of pairwise isospectral maps. By the second property and (13), the function $t \mapsto \|\mathrm{Ric}^{G_{j(t)}}|_{\mathfrak{v}}\|^2$

(where $j(t) := F_Y^t(j)$) has nonzero derivative at $t = 0$ for every $j$ in the Zariski open subset $U := \mathcal{J} \setminus \{j \in \mathcal{J} \mid dq|_j(Y) = 0\} \neq \emptyset$. The last statement of the proposition now follows by Proposition 3.5 and the fact that $\|\mathrm{Ric}^{G_j}|_{\mathfrak{v}}\|^2$ equals $\lambda_1^2 + \ldots \lambda_m^2$, where $\lambda_1, \ldots, \lambda_m$ are the eigenvalues of $\mathrm{Ric}^{G_j}|_{\mathfrak{v}}$. □

4.2 LEMMA. *Let $p_{a,b}$ and $q$ be the polynomials on the space $\mathcal{J}$ of all linear maps from $\mathbb{R}^2$ to $\mathrm{so}(m)$ which were defined in* (11) *and* (14). *Then the vector field $Y$ on $\mathcal{J}$, given by $Y(j) = (j_1^3 j_2 - j_2 j_1^3, 0)$, satisfies $dp_{a,b}(Y) = 0$ for all $a, b \in \mathbb{N}_0$ such that $a + b > 0$, while the polynomial $j \to dq|_j(Y)$ does not vanish identically on $\mathcal{J}$ if $m \geq 5$.*

*Proof.* First note that for all $a, b \in \mathbb{N}$, we have $dp_{a,0}(Y) = a \cdot \mathrm{tr}\big(j_1^{a-1}(j_1^3 j_2 - j_2 j_1^3)\big) = 0$, and trivially $dp_{0,b}(Y) = 0$. Thus we can assume in the following that both $a$ and $b$ are nonzero. Let $j \in \mathcal{J}$, and let $\varepsilon = (\varepsilon_1, \varepsilon_2) \in T_j \mathcal{J} \cong \mathcal{J}$. Then

$$dp_{a,b}|_j(\varepsilon) = \sum_{\sigma \in \mathfrak{S}_{a,b}} \sum_{i=1}^{a+b} \mathrm{tr}\big(j_{\sigma(1)} \cdots j_{\sigma(i-1)} \varepsilon_{\sigma(i)} j_{\sigma(i+1)} \cdots j_{\sigma(a+b)}\big)$$

$$= \sum_{i=1}^{a+b} \sum_{\sigma \in \mathfrak{S}_{a,b}} \mathrm{tr}\big(\varepsilon_{\sigma(i)} j_{\sigma(i+1)} \cdots j_{\sigma(a+b)} j_{\sigma(1)} \cdots j_{\sigma(i-1)}\big)$$

$$= (a+b) \cdot \Big( \sum_{\tau \in \mathfrak{S}_{a-1,b}} \mathrm{tr}(\varepsilon_1 j_{\tau(1)} \cdots j_{\tau(a+b-1)})$$

$$+ \sum_{\tau \in \mathfrak{S}_{a,b-1}} \mathrm{tr}(\varepsilon_2 j_{\tau(1)} \cdots j_{\tau(a+b-1)}) \Big).$$

Showing $dp_{a,b}(Y) = 0$ thus reduces to showing that

(15) $$\sum_{\tau \in \mathfrak{S}_{a-1,b}} \mathrm{tr}(j_1^3 j_2 j_{\tau(1)} \cdots j_{\tau(a+b-1)}) = \sum_{\tau \in \mathfrak{S}_{a-1,b}} \mathrm{tr}(j_2 j_1^3 j_{\tau(1)} \cdots j_{\tau(a+b-1)})$$

for all $j \in \mathcal{J}$. We now define a bijection $\varphi$ from $\mathfrak{S}_{a-1,b}$ to itself such that for $\rho = \varphi(\tau)$,

$$\mathrm{tr}(j_1^3 j_2 j_{\tau(1)} \cdots j_{\tau(a+b-1)}) = \mathrm{tr}(j_2 j_1^3 j_{\rho(1)} \cdots j_{\rho(a+b-1)}).$$

Namely, given $\tau \in \mathfrak{S}_{a-1,b}$, let $c = \max\{i \mid \tau(i) = 2\}$; define $\rho = \varphi(\tau)$ by letting $\rho(i) = 1$ for $i = 1, \ldots, a+b-c-1$, $\rho(a+b-c) = 2$, and $\rho(a+b-c+i) = \tau(i)$ for $i = 1, \ldots, c-1$. Obviously $\varphi : \mathfrak{S}_{a-1,b} \to \mathfrak{S}_{a-1,b}$ is a bijection; moreover,

$$\mathrm{tr}(j_1^3 j_2 j_{\tau(1)} \cdots j_{\tau(a+b-1)}) = \mathrm{tr}(j_1^3 j_2 j_{\tau(1)} \cdots j_{\tau(c-1)} j_2 j_1^{a+b-c-1})$$

$$= \mathrm{tr}(j_2 j_1^3 j_1^{a+b-c-1} j_2 j_{\tau(1)} \cdots j_{\tau(c-1)})$$

$$= \mathrm{tr}(j_2 j_1^3 j_{\rho(1)} \cdots j_{\rho(a+b-1)}),$$

as required. Equation (15) now follows; hence $dp_{a,b}(Y) = 0$ as claimed. Moreover,

$$dq|_j Y = \text{tr}\big((j_1^3 j_2 - j_2 j_1^3) j_1 j_2^2 + j_1(j_1^3 j_2 - j_2 j_1^3) j_2^2\big) = \text{tr}(j_1^3 j_2 j_1 j_2^2 - j_1 j_2 j_1^3 j_2^2).$$

For $m \geq 5$, this polynomial does not vanish identically on $\mathcal{J}$; e.g., for

$$(16) \qquad j_1 = \begin{pmatrix} 0 & 0 & 0 & 1 & 0 \\ 0 & 0 & 0 & 1 & 1 \\ 0 & 0 & 0 & 0 & 0 \\ -1 & -1 & 0 & 0 & 0 \\ 0 & -1 & 0 & 0 & 0 \end{pmatrix} \text{ and } j_2 = \begin{pmatrix} 0 & 1 & 0 & 0 & 0 \\ -1 & 0 & 0 & 0 & 0 \\ 0 & 0 & 0 & 1 & 0 \\ 0 & 0 & -1 & 0 & 0 \\ 0 & 0 & 0 & 0 & 0 \end{pmatrix}$$

it equals $2 \neq 0$. $\square$

4.3 *Remarks.* (i) In case $m \leq 4$, we have $\text{tr}(j_1^3 j_2 j_1 j_2^2 - j_1 j_2 j_1^3 j_2^2) = 0$ for all $j_1, j_2$ in $so(m)$, which is not hard to show. We also remark that on the other hand, this is implied by the above arguments in connection with Remark 2.4(iii).

(ii) In Example 2.5, it happens to be the case that $\|\text{Ric}^{G_{j(t)}}|_{\mathfrak{v}}\|^2$ is constant in $t$. In other words, even though $\det(\text{Ric}^{G_{j(t)}}|_{\mathfrak{v}})$ is nonconstant in $t$ (as observed in 2.5) and thus $\text{Ric}^{G_{j(t)}}$ has nonconstant eigenvalue spectrum, not only the trace of $\text{Ric}^{G_{j(t)}}|_{\mathfrak{v}}$ is constant in $t$ (which is automatic by (10) since the $j(t)$ are pairwise isospectral), but even the sum of the squares of its eigenvalues is constant. Hence the families given in 2.5 do not represent examples for Proposition 4.1. Examples for 4.1 would be obtained by computing flow lines of the vector field $Y$ through suitable starting points $j = (j_1, j_2)$, e.g., the one given by (16). However, this seems hard to do explicitly. A different explicit example of an isospectral family $j(t) = (j_1(t), j_2(t)) : \mathbb{R}^2 \to so(5)$ with $\|\text{Ric}^{G_{j(t)}}|_{\mathfrak{v}}\|^2 \neq \text{const}$, although not representing a flow line of $Y$, is given by

$$j_1(t) := \begin{pmatrix} 0 & 0 & -t & 0 & 0 \\ 0 & 0 & 0 & t-1 & 0 \\ t & 0 & 0 & 0 & -\varphi(t) \\ 0 & 1-t & 0 & 0 & -\psi(t) \\ 0 & 0 & \varphi(t) & \psi(t) & 0 \end{pmatrix}, \qquad j_2(t) := \begin{pmatrix} 0 & 1 & 0 & 0 & 0 \\ -1 & 0 & 0 & 0 & 0 \\ 0 & 0 & 0 & 1 & 0 \\ 0 & 0 & -1 & 0 & 0 \\ 0 & 0 & 0 & 0 & 0 \end{pmatrix},$$

where $\varphi(t) = \big((t^4 - 3t^2 + 1)/(1 - 2t)\big)^{1/2}$ and $\psi(t) = \big((-t^4 + 4t^3 - 3t^2 - 2t + 1)/(1-2t)\big)^{1/2}$. This deformation is defined for $t \in \big[\frac{1}{2}(1-\sqrt{5}), \frac{1}{2}(3-\sqrt{5})\big]$. The $j(t)$ are pairwise isospectral since $\det\big(\lambda \text{Id} - (sj_1(t) + uj_2(t))\big) = \lambda^5 + (3s^2 + 2u^2)\lambda^3 + (s^2 + u^2)\lambda$ is independent of $t$. On the other hand, $\|\text{Ric}^{G_{j(t)}}|_{\mathfrak{v}}\|^2 = t^2 - t + \frac{13}{2}$ is nonconstant in $t$.

## 5. The spectrum on 1-forms

Let $(M, g)$ be a compact Riemannian manifold without boundary. The Hodge-de Rham Laplacian $\Delta_g^p$ acts on the space of smooth $p$-forms on $(M, g)$ by $\Delta_g^p(\alpha) = -(d\delta + \delta d)\alpha$, where $d$ denotes the exterior derivative and $\delta$ its adjoint operator with respect to the canonical scalar product on $p$-forms induced by $g$. Like $\Delta_g = \Delta_g^0$, the operators $\Delta_g^p$ are elliptic, positive semi-definite and self-adjoint, and have a discrete spectrum of eigenalues tending to infinity. Two manifolds $(M, g)$ and $(M', g')$ are called *isospectral on $p$-forms* if $\operatorname{spec}(\Delta_g^p) = \operatorname{spec}(\Delta_{g'}^p)$.

5.1 PROPOSITION. *If $j \sim j'$ and $\|\operatorname{Ric}^{G_j}|_{\mathfrak{v}}\|^2 \neq \|\operatorname{Ric}^{G_{j'}}|_{\mathfrak{v}}\|^2$ then the isospectral manifolds $(C, g_j)$ and $(C, g_{j'})$ are not isospectral on 1-forms. In particular, in each of the families $(C, g_{j(t)})$ from Proposition 4.1, the manifolds are not pairwise isospectral on 1-forms.*

We will prove this proposition by using the heat invariants for the spectrum on functions and on 1-forms; see [Gi, §4.8] for details. The heat invariants $a_i^p(g)$ are the coefficients appearing in the asymptotic expansion

$$\operatorname{tr}\bigl(\exp(-t\Delta_g^p)\bigr) \;\sim\; (4\pi t)^{-\dim M/2} \sum_{i=0}^{\infty} a_i^p(g) t^i \quad \text{for } t \searrow 0.$$

If $(M, g)$ and $(M', g')$ are isospectral on $p$-forms, then necessarily $a_i^p(g) = a_i^p(g')$ for all $i \in \mathbb{N}_0$. We will use this fact for the heat invariants $a_2^0(g)$ and $a_2^1(g)$ which are given explicitly as follows.

5.2 LEMMA ([Gi, Theorem 4.8.18]). *Let $(M, g)$ be a closed Riemannian manifold, and let $c_s(g) := \int_M \operatorname{scal}^2(p)\, d\mathrm{vol}_g(p)$, $c_{\operatorname{Ric}}(g) := \int_M \|\operatorname{Ric}\|_p^2\, d\mathrm{vol}_g(p)$, and $c_R(g) := \int_M \|R\|_p^2\, d\mathrm{vol}_g(p)$, where $\operatorname{scal}$ denotes the scalar curvature, $\operatorname{Ric}$ and $R$ denote the Ricci and curvature tensors of $(M, g)$, respectively, and the norm $\|.\|_p$ is the euclidean norm on tensors, corresponding to the scalar product induced by $g$ on $T_pM$. Then*

(i) $a_2^0(g) = \frac{1}{360}(5c_s - 2c_{\operatorname{Ric}} + 2c_R)(g)$,

(ii) $a_2^1(g) = \frac{1}{360}\bigl((5\dim M - 60)c_s - (2\dim M - 180)c_{\operatorname{Ric}} + (2\dim M - 30)c_R\bigr)(g)$.

5.3 COROLLARY. *If $(M, g)$ and $(M', g')$ are isospectral on functions and if $a_2^1(g) = a_2^1(g')$ (in particular, if $(M, g)$ and $(M', g')$ are isospectral both on functions and on 1-forms), then*

$$c_s(g) - c_s(g') + 10\bigl(c_{\operatorname{Ric}}(g) - c_{\operatorname{Ric}}(g')\bigr) = 0.$$

*Proof.* Since $a_2^0(g) - a_2^0(g') = 0$ by isospectrality on functions, we have by Lemma 5.2 that $c_R(g) - c_R(g') = -\frac{5}{2}(c_s(g) - c_s(g')) + c_{\operatorname{Ric}}(g) - c_{\operatorname{Ric}}(g')$. Substituting this into $a_2^1(g) - a_2^1(g') = 0$ yields the required formula. □

The next lemma is crucial for the proof of Proposition 5.1.

5.4 LEMMA. *If $j \sim j'$ then for the metrics $g_j$ and $g_{j'}$ on $C = S^1(\mathfrak{v}) \times S$ the following holds:*

(i) $c_s(g_j) - c_s(g_{j'}) = \text{vol}(S) \cdot \text{vol}(S^{m-1}) \cdot \dfrac{2}{m(m+2)} \cdot \left(\|\text{Ric}^{G_j}|_{\mathfrak{v}}\|^2 - \|\text{Ric}^{G_{j'}}|_{\mathfrak{v}}\|^2\right),$

(ii) $c_{\text{Ric}}(g_j) - c_{\text{Ric}}(g_{j'}) = \text{vol}(S) \cdot \text{vol}(S^{m-1}) \cdot \dfrac{m-4}{m} \cdot \left(\|\text{Ric}^{G_j}|_{\mathfrak{v}}\|^2 - \|\text{Ric}^{G_{j'}}|_{\mathfrak{v}}\|^2\right).$

Here $m = \dim \mathfrak{v}$, and $\text{vol}(S)$ is the volume of $(S, k)$, where $k$ is the biinvariant metric used in 1.3(ii).

*Proof of Proposition* 5.1. By Lemma 5.4,

$$c_s(g_j) - c_s(g_{j'}) + 10\big(c_{\text{Ric}}(g_j) - c_{\text{Ric}}(g_{j'})\big)$$
$$= \text{vol}(S) \cdot \text{vol}(S^{m-1}) \cdot \dfrac{10m^2 - 20m - 78}{m(m+2)} \cdot \big(\|\text{Ric}^{G_j}|_{\mathfrak{v}}\|^2 - \|\text{Ric}^{G_{j'}}|_{\mathfrak{v}}\|^2\big).$$

This is nonzero since the third factor is nonzero for all integers $m$ and the last factor is nonzero by assumption. Since $(C, g_j)$ and $(C, g_{j'})$ are isospectral on functions by Theorem 2.2, Corollary 5.3 implies that $a_2^1(g_j) \neq a_2^1(g_{j'})$; hence the two manifolds are not isospectral on 1-forms. $\square$

*Proof of Lemma* 5.4. (i) By Lemma 3.3,

$$\int_C (\text{scal}_p^{C,j})^2 \, d\text{vol}_{g_j}(p)$$
$$= \text{vol}(S) \cdot \int_{S^1(\mathfrak{v})} \big(\text{Ric}^{G_j}(v,v)^2 - 2\varphi(j)\text{Ric}^{G_j}(v,v) + \varphi(j)^2\big) \, dv,$$

where $\varphi(j) = \text{scal}^{G_j} + \text{scal}^S + (m-1)(m-2)$ does not depend on $v$. Recall from the proof of 3.5 that $j \sim j'$ implies $\varphi(j) = \varphi(j')$. Moreover, for every $A \in \text{End}(\mathbb{R}^m)$ we have the general formulas

(17)
$$\int_{S^{m-1}} \langle Av, v \rangle \, dv = \text{tr}(A) \cdot \int_{S^{m-1}} v_1^2 \, dv = \text{tr}(A) \cdot \text{vol}(S^{m-1}) \cdot \dfrac{1}{m},$$
$$\int_{S^{m-1}} \langle Av, v \rangle^2 \, dv = \big(\text{tr}(A)^2 + \|A\|^2 + \langle A, {}^tA\rangle\big) \cdot \int_{S^{m-1}} v_1^2 v_2^2 \, dv$$
$$= \big(\text{tr}(A)^2 + \|A\|^2 + \langle A, {}^tA\rangle\big) \cdot \text{vol}(S^{m-1}) \cdot \dfrac{1}{m(m+2)},$$

where $v_1, \ldots, v_m$ are the standard coordinates on $\mathbb{R}^m$. If $j \sim j'$ then $\text{tr}(\text{Ric}^{G_j}|_{\mathfrak{v}})$ $= \frac{1}{2} \sum_{i=1}^r \text{tr}(j_{z_i}^2)$ equals $\text{tr}(\text{Ric}^{G_{j'}}|_{\mathfrak{v}}) = \frac{1}{2} \sum_{i=1}^r \text{tr}(j'^2_{z_i})$. Since $\text{Ric}^{G_j}|_{\mathfrak{v}} : \mathfrak{v} \to \mathfrak{v}$ is symmetric, (17) therefore implies that $c_s(g_j) - c_s(g_{j'}) = \text{vol}(S) \cdot \text{vol}(S^{m-1}) \cdot \dfrac{1}{m(m+2)} \cdot \big(2\|\text{Ric}^{G_j}|_{\mathfrak{v}}\|^2 - 2\|\text{Ric}^{G_{j'}}|_{\mathfrak{v}}\|^2\big)$, as claimed.

(ii) In order to access $c_{\text{Ric}}(g_j)$ we first establish explicit formulas for the Ricci tensor of $(C, g_j)$, denoted $\text{Ric}^{C,j}$. Recall from the proof of 3.3 that for any $(v, s) \in S^1(\mathfrak{v}) \times S = C$, the tangent space $T_{(v,s)}C$ is given as $\{(x^j + h^j + u^j)_{(v,s)} \mid x \in \mathfrak{v}, x \perp v, h \in \mathfrak{h}, u \in \mathfrak{u} = \mathfrak{h}^\perp \cap \mathfrak{s}\}$, where $x^j$ etc. are the vector fields defined in 1.7. Using Lemmas 3.1 and 1.6, the fact that $R(u^j, v^j)v^j = 0$ in $(E, g_j)$ by the proof of 3.2, the results on the shape operator $\nabla N$ from the proof of 3.3, formula (7), and the Gauss equations, one easily shows for all $x, x' \in \mathfrak{v}$ such that $x, x' \perp v$, all $u, u' \in \mathfrak{u}$, all $h, h' \in \mathfrak{h}$ and the corresponding vectors $z = f^{-1}(h)$, $z' = f^{-1}(h') \in \mathfrak{z}$:

$$\begin{aligned}
\mathrm{Ric}^{C,j}_{(v,s)}(x^j, x'^j) &= \langle \mathrm{Ric}^{G_j} x, x' \rangle + \frac{1}{2}\langle [x, v], [x', v] \rangle + (m-2)\langle x, x' \rangle, \\
\mathrm{Ric}^{C,j}_{(v,s)}(x^j, h^j) &= \frac{m-2}{2}\langle j_z x, v \rangle, \\
\mathrm{Ric}^{C,j}_{(v,s)}(x^j, u^j) &= 0, \\
\mathrm{Ric}^{C,j}_{(v,s)}(h^j, h'^j) &= -\frac{1}{4}\mathrm{tr}(j_z j_{z'}) - \frac{1}{2}\langle j_z v, j_{z'} v \rangle + \mathrm{Ric}^S(h, h'), \\
\mathrm{Ric}^{C,j}_{(v,s)}(h^j, u^j) &= \mathrm{Ric}^S(h, u), \\
\mathrm{Ric}^{C,j}_{(v,s)}(u^j, u'^j) &= \mathrm{Ric}^S(u, u'),
\end{aligned}$$

where $\mathrm{Ric}^S$ denotes the Ricci tensor of $(S, k)$. Thus

(18)
$$\begin{aligned}
\|\mathrm{Ric}^{C,j}_{(v,s)}\|^2 &= \sum_{a,b=1}^m \Big(\langle \mathrm{Ric}^{G_j} x_a, x_b \rangle + \frac{1}{2}\langle [x_a, v], [x_b, v] \rangle + (m-2)\langle x_a, x_b \rangle\Big)^2 \\
&\quad - 2\sum_{a=1}^m \big(\langle \mathrm{Ric}^{G_j} x_a, v \rangle + (m-2)\langle x_a, v \rangle\big)^2 \\
&\quad + \big(\langle \mathrm{Ric}^{G_j} v, v \rangle + (m-2)\big)^2 \\
&\quad + 2\sum_{a=1}^m \sum_{i=1}^r \Big(\frac{m-2}{2}\langle j_{z_i} x_a, v \rangle^2\Big) \\
&\quad + \sum_{i,k=1}^r \Big(-\frac{1}{4}\mathrm{tr}(j_{z_i} j_{z_k}) - \frac{1}{2}\langle j_{z_i} v, j_{z_k} v \rangle\Big)^2 \\
&\quad + 2\sum_{i,k=1}^r \Big(-\frac{1}{4}\mathrm{tr}(j_{z_i} j_{z_k}) - \frac{1}{2}\langle j_{z_i} v, j_{z_k} v \rangle\Big) \cdot \mathrm{Ric}^S(h_i, h_k) \\
&\quad + \|\mathrm{Ric}^S\|^2,
\end{aligned}$$

where $\{x_1, \ldots, x_m\}$ and $\{z_1, \ldots, z_r\}$ are orthonormal bases of $\mathfrak{v}$ and $\mathfrak{z}$, respectively, and $h_i = f(z_i)$ for $i = 1, \ldots, r$. By $\|\mathrm{Ric}^S\|$ we denote the pointwise norm of $\mathrm{Ric}^S$ which is constant since $(S, k)$ is homogeneous; moreover, this norm does not depend on $j$. Note also that $j \sim j'$ implies $\mathrm{tr}(j_z j_w) = \mathrm{tr}(j'_z j'_w)$ as well as $\mathrm{tr}(2 j_z^2 j_w^2 + j_z j_w j_z j_w) = \mathrm{tr}(2 j'^2_z j'^2_w + j'_z j'_w j'_z j'_w)$ for all $z, w \in \mathfrak{z}$; this follows by considering the coefficients at $su$, resp. $s^2 u^2$, in $\mathrm{tr}\big((sj_z + uj_w)^2\big)$, resp. $\operatorname{tr}((sj_z + uj_w)^4)$. Using these facts and formulas (10) and (17), one derives from (18) by direct computation (of which we omit the details here) that the only summands in $\int_C \|\operatorname{Ric}^{C,j}\|_p^2 \, d\operatorname{vol}_{g_j}(p)$ which in case $j \sim j'$ do possibly *not* coincide for $j$ and $j'$ consist of linear multiples of $\|\operatorname{Ric}^{G_j}|_{\mathfrak{v}}\|^2$; moreover, the respective coefficients sum up to $\operatorname{vol}(S) \cdot \int_{S^1(\mathfrak{v})} (1 - 4v_1^2 + 0\, v_1^2 v_2^2)\, dv = \operatorname{vol}(S) \cdot \operatorname{vol}(S^{m-1}) \cdot (1 - \frac{4}{m})$. □

5.5 *Final Remarks.* (i) The proof of Proposition 5.1 shows that the heat invariant $a_2^1(g_{j(t)})$ changes nontrivially during each of the isospectral deformations given in Proposition 4.1. As mentioned in the introduction, these families consitute the first examples of manifolds which are isospectral on functions and for which nonisospectrality on 1-forms turns out already from considering heat invariants only. Not also that the terms $c_s(g_{j(t)})$, $c_{\operatorname{Ric}}(g_{j(t)})$, and $c_R(g_{j(t)})$ change nontrivially as $t$ varies. The heat invariant $a_2^0(g_{j(t)})$ which is a linear combination of these (see 5.2(i)) is nevertheless constant in $t$ by isospectrality on functions, while $a_2^1(g_{j(t)})$ (another linear combination of the same terms, see 5.2(ii)) is not constant in $t$.

(ii) Note that Proposition 5.1 does not exclude the possible existence of isospectral manifolds $(C, g_j)$ and $(C, g_{j'})$ with $j \sim j'$ which are also isospectral on 1-forms. We do not know whether such examples exist here. The proposition only says that in this case, the norms of $\operatorname{Ric}^{G_j}|_{\mathfrak{v}}$ and $\operatorname{Ric}^{G_{j'}}|_{\mathfrak{v}}$ must coincide, which by 4.1 is frequently not the case.

(iii) Each of the conditions $j \sim j'$ and $\|\operatorname{Ric}^{G_j}|_{\mathfrak{v}}\|^2 \ne \|\operatorname{Ric}^{G_{j'}}|_{\mathfrak{v}}\|^2$ on a pair $j, j'$ of linear maps from $\mathfrak{z}$ to $\operatorname{so}(\mathfrak{v})$ is invariant under rescaling both maps by the same real factor. Lemma 1.7 implies that the metrics $g_{\alpha j}$ converge to the product metric on $C = S^1(\mathfrak{v}) \times S$ for $\alpha \to 0$, where $S^1(\mathfrak{v})$ is the round unit sphere in $\mathfrak{v} = \mathbb{R}^m$, and $S$ is endowed with the biinvariant metric $k$ used in 1.3(ii). Therefore, it follows from Proposition 4.1 that there are nontrivial isospectral families of metrics $g_{j(t)}$ on $S^n \times S$ for every $n \ge 4$ which are arbitrarily close to the product metric $g_0$. In particular, nontrivial isospectral deformations of the form $(S^n \times S^3 \times S^3, g_{j(t)})$ can take place arbitrarily close to the standard metric on $S^n \times S^3 \times S^3$ for every $n \ge 4$.

MATHEMATISCHES INSTITUT, UNIVERSITÄT BONN, GERMANY
*E-mail address*: schueth@math.uni-bonn.de